\documentclass[11pt]{amsart}
\usepackage[T1]{fontenc}
\usepackage{lmodern}
\usepackage{amsmath,amssymb,amsthm}
\usepackage[all]{xy}
\usepackage{microtype}
\usepackage[margin=1.05in]{geometry}
\usepackage[hidelinks]{hyperref}
\newtheorem{theorem}{Theorem}[section]
\newtheorem{proposition}[theorem]{Proposition}
\newtheorem{lemma}[theorem]{Lemma}
\newtheorem{corollary}[theorem]{Corollary}
\theoremstyle{definition}
\newtheorem{definition}[theorem]{Definition}
\newtheorem{example}[theorem]{Example}
\theoremstyle{remark}
\newtheorem{remark}[theorem]{Remark}
\DeclareMathOperator{\Cl}{Cl}
\DeclareMathOperator{\Pic}{Pic}
\DeclareMathOperator{\Spec}{Spec}
\DeclareMathOperator{\Ann}{Ann}
\DeclareMathOperator{\Supp}{Supp}
\DeclareMathOperator{\Sat}{Sat}
\DeclareMathOperator{\qf}{qf}
\newcommand{\NN}{\mathbb N}
\newcommand{\ZZ}{\mathbb Z}
\newcommand{\QQ}{\mathbb Q}

\newcommand{\tMax}{t\text{-}\operatorname{Max}}
\newcommand{\doi}[1]{\href{https://doi.org/#1}{\nolinkurl{#1}}}

\title[Almost $S$-GCD domains]{Almost $S$-GCD domains: localization, class groups, and extensions}
\author[A. Hamed]{Ahmed Hamed}
\address{(Hamed) Department of Mathematics, Faculty of Sciences, Monastir, Tunisia}
\email{hamed.ahmed@hotmail.fr}
\author[H. Kim]{Hwankoo Kim}
\address{(Kim) Division of Computer Engineering, Hoseo University, Asan 31499, Republic of Korea}
\email{hkkim@hoseo.edu}
\author[M.T. Ahmed]{Malik Tusif Ahmed}
\address{(Ahmed) Department of Mathematics \& Statistics, Hazara University,
Dhodial, Mansehra, Khyber Pakhtunkhwa 21120, Pakistan}
\email{maliktusifahmed@outlook.com}
\thanks{Corresponding author: H. Kim.}
\date{\today}

\begin{document}
\begin{abstract}
Let $D$ be an integral domain and let $S\subseteq D\setminus\{0\}$ be
multiplicative. We study the condition that, for each nonzero $a,b\in D$,
some $(a^n,b^n)_v$ is $S$-principal. A power estimate extends this binary
condition to arbitrary finite families and permits synchronization of
exponents. We describe the obstruction to descent from the AGCD domain
$D_S$ by bounded $S$-torsion and by a finiteness condition on intersections
of principal ideals. For an almost Pr\"ufer $v$-multiplication domain,
the almost $S$-GCD property, the AGCD property of $D_S$, and torsion of the
$S$-class group are equivalent. In Krull domains this gives a criterion
in terms of the divisor classes of height-one primes meeting $S$,
and finite rational class rank reduces the multiplicative set to one
element. We prove conductor transfer without a principal-contraction
hypothesis, polynomial transfer in several settings, a normalization
homeomorphism away from $S$, and a finite-principal-cover criterion with
a Picard-group obstruction. Explicit examples distinguish localization,
bounded descent, and the global and sublocal conditions.
\end{abstract}
\maketitle

\section{Introduction}\label{intro}

The theory of almost GCD domains replaces the existence of a greatest
common divisor for a pair by its existence after raising both elements
to a suitable common power. This connects divisibility with torsion in
ideal class groups and with root extensions to integral closures.
Zafrullah introduced AGCD domains in \cite{Z3}; Anderson and Zafrullah
developed their relation to almost B\'ezout domains in \cite{AZ-Bez}.
Almost splitting sets and the additional obstruction to polynomial
stability were studied in \cite{ATZ,AndersonChang2007}. The distinction
between local principality and global principality is developed further
in the theory of locally and almost generalized GCD domains
\cite{ACZlocal}.

A different relaxation allows a multiplicative set to absorb the
failure of an ideal to be principal. A nonzero fractional ideal $I$ is
\emph{$S$-principal} if $sI\subseteq dD\subseteq I$ for some $s\in S$
and $d\in I\setminus\{0\}$. Hamed and Hizem introduced $S$-GCD domains
and the $S$-class group in \cite{Hamed}; Anderson, Hamed, and Zafrullah
gave a systematic study of $S$-GCD domains in \cite{AHZ}.
These constructions build on the $S$-finiteness viewpoint of
\cite{Anderson1} and the class-group framework of \cite{Alain}.

We combine the power condition with $S$-principality. The resulting
notion is sensitive to a uniform denominator: principality after
localization only clears denominators element by element, whereas
$S$-principality requires one element of $S$ to work on an entire ideal.
Consequently, descent from an AGCD localization can fail even when
the original domain is integrally closed. Understanding when this
obstruction disappears is the main organizing question of this paper.

Our first contribution is an explicit power estimate that proves the
finite-family form of the definition without normality or finiteness
hypotheses (Theorem~\ref{t:finite-families}). We then give two complementary
descent mechanisms: bounded torsion in principal quotients, and
$S$-$v$-finiteness of intersections of principal ideals
(Theorems~\ref{t3} and \ref{t:finite-conductor-localization}).
The latter separates the actual finiteness obstruction from the
stronger contraction condition.

The class-group results identify a setting in which no extra
contraction assumption is needed. For an almost Pr\"ufer
$v$-multiplication domain $D$, we prove the equivalence of the
almost $S$-GCD property, the AGCD property of $D_S$, and torsion of
$S\text{-}\Cl_t(D)$ (Theorem~\ref{t:apvmd-relative-characterization}).
The proof includes the localization map for $t$-invertible ideals,
with its precise kernel, in the context of the established relative
class-group theory \cite{HamedNagata, Hamed}.
For Krull domains, the classical localization sequence for divisors
\cite{Fossum} yields a criterion using the classes of primes meeting
$S$, together with a single-denominator consequence when
$\Cl(D)\otimes_{\ZZ}\QQ$ is finite dimensional.

We also study extensions and local constructions. Over an integrally
closed domain the principal-contraction condition is equivalent to
the corresponding condition in any polynomial ring
(Proposition~\ref{p:normal-polynomial-contraction}).
For PvMDs, polynomial transfer holds for an arbitrary set of
indeterminates and an arbitrary $S$.
If an overring has a conductor meeting $S$, the almost $S$-GCD
property transfers in both directions
(Theorem~\ref{t:conductor-transfer}). This gives pullback and
normalization applications. Finally, finite principal covers produce
invertible power intersections, and the Picard group measures the
remaining obstruction to their principality
(Theorem~\ref{t:principal-cover-picard}).
The results use the classical AGCD, splitting-set, and
divisor-localization theorems as background.

\subsection{Conventions and preliminaries}

All rings are commutative with identity, all domains are nonzero, and
all inclusions of rings preserve identity. For a domain $D$, put
$D^\ast=D\setminus\{0\}$ and $K=\qf(D)$. Every multiplicative set
contains $1$ and is contained in $D^\ast$. Write $U(D)$ for the
group of units, $\NN=\{1,2,\ldots\}$, and
$\langle s\rangle=\{s^j:j\geq0\}$; the localization at this cyclic
set is denoted by $D_s$.

For a nonzero fractional ideal $I$ of $D$, set
$I^{-1}=(D:I)=\{x\in K:xI\subseteq D\}$ and $I_v=(I^{-1})^{-1}$.
The associated finite-type operation is
\[
 I_t=\bigcup\{F_v:0\ne F\subseteq I,\ F\text{ finitely generated}\}.
\]
An ideal is \emph{divisorial} if $I=I_v$, and is a \emph{$t$-ideal}
if $I=I_t$. A divisorial ideal is \emph{of $v$-finite type} if it is
$F_v$ for a finitely generated fractional ideal $F$. When $F$ is
finitely generated, $F_t=F_v$. We use subscripts for all star operations;
an operation on an overring is marked explicitly when necessary.
Standard properties are recalled from \cite{Alain, Gilmer}.

A fractional ideal $I$ is \emph{$t$-invertible} if $(II^{-1})_t=D$.
A $t$-invertible $t$-ideal and its inverse are divisorial and of
$v$-finite type. Maximal proper $t$-ideals are prime; their set is
denoted by $\tMax(D)$. Every fractional $t$-ideal is the intersection
of its extensions to $D_M$ for $M\in\tMax(D)$.
A \emph{Pr\"ufer $v$-multiplication domain} (PvMD) is a domain in which
every nonzero finitely generated ideal is $t$-invertible.
Equivalently, $D_M$ is a valuation domain for every $M\in\tMax(D)$.
An \emph{almost Pr\"ufer $v$-multiplication domain} (APvMD) is a domain
in which, for each $a,b\in D^\ast$, some $(a^n,b^n)$ is $t$-invertible
\cite{LiAPvMD}. Every PvMD is an APvMD.

The domain $D$ is \emph{AGCD} if, for each $a,b\in D^\ast$, some
$a^nD\cap b^nD$ is principal. Equivalently, some $(a^n,b^n)_v$ is
principal. It is an \emph{$S$-GCD domain} if every nonzero finitely
generated ideal $F$ has $S$-principal $v$-closure. Scalar multiplication
by an element of $K^\ast$ preserves $S$-principality.
We write $\overline D$ for the integral closure of $D$ in $K$.
An extension $A\subseteq B$ is a \emph{root extension} if every
$b\in B$ has a positive power in $A$.

\section{Finite families and descent from localization}\label{sec:basic}

We first justify passage from pairs to finite families and synchronize
the exponents. We then isolate the difference between localized
principality and $S$-principality. This yields descent criteria and
examples showing which hypotheses cannot simply be dropped.

\subsection{Powers, intersections, and finite families}

\begin{definition}\label{d:almost-S}
Let $D$ be a domain and let $S\subseteq D^\ast$ be multiplicative.
The domain $D$ is \emph{almost $S$-GCD} if, for each $a,b\in D^\ast$,
there exists $n\geq1$ for which $(a^n,b^n)_v$ is $S$-principal.
\end{definition}

Every $S$-GCD domain is almost $S$-GCD. If $S\subseteq U(D)$, the
definition is the ordinary AGCD condition. At the other extreme,
every domain is an $S$-GCD domain for $S=D^\ast$: for a nonzero
integral ideal $I$ and $d\in I\setminus\{0\}$, one has
$dI\subseteq dD\subseteq I$.

The basic two-element characterization is the power analogue of
\cite[Theorem~2.1]{AHZ}.

\begin{theorem}\label{t1}
Let $D$ be a domain and $S\subseteq D^\ast$ multiplicative. The
following conditions are equivalent.
\begin{enumerate}
\item $D$ is almost $S$-GCD.
\item For every $a,b\in D^\ast$, some $a^nD\cap b^nD$ is $S$-principal.
\item For every $a,b\in D^\ast$, some $(a^nD:b^nD)$ is $S$-principal.
\end{enumerate}
For a fixed pair, the same exponents can be used in the three conditions.
\end{theorem}

\begin{proof}
If $I$ is a nonzero fractional divisorial ideal and
$sI\subseteq dD\subseteq I$, inversion gives
$sI^{-1}\subseteq(s/d)D\subseteq I^{-1}$.
Since $(I^{-1})^{-1}=I$, inversion preserves $S$-principality in
both directions. Now
\[
 (a^n,b^n)^{-1}=(a^nb^n)^{-1}(a^nD\cap b^nD),
 \qquad
 a^nD\cap b^nD=b^n(a^nD:b^nD).
\]
Scalar multiplication and inversion therefore prove the assertions.
\end{proof}

The next estimate is needed when different pairs initially yield
different exponents.

\begin{lemma}\label{l:power-sandwich}
Let $D$ be a domain, let $S\subseteq D^\ast$ be multiplicative, and let
$a_1,\ldots,a_r\in D^\ast$. Suppose that
$s(a_1,\ldots,a_r)_v\subseteq dD\subseteq(a_1,\ldots,a_r)_v$
for some $s\in S$ and $d\in D^\ast$. For every $m\geq1$, put
$q=(r-1)(m-1)$. Then
\[
s^q d^mD\subseteq(a_1^m,\ldots,a_r^m)_v
\subseteq s^{-m}d^mD.
\]
In particular, $(a_1^m,\ldots,a_r^m)_v$ is $S$-principal.
\end{lemma}

\begin{proof}
Set $u_i=sa_i/d\in D$ and $U=(u_1,\ldots,u_r)D$.
The given inclusions imply $s\in U_v$ and $U\subseteq D$.
For $N=r(m-1)+1=m+q$, the pigeonhole principle gives
$U^N\subseteq(u_1^m,\ldots,u_r^m)D$. Since
$((U_v)^N)_v=(U^N)_v$, we obtain
$s^N\in(u_1^m,\ldots,u_r^m)_v$.
Multiplying by $d^m/s^m$ proves the left inclusion; the right one follows
from $u_i^m\in D$. Consequently
\[
s^{m+q}(a_1^m,\ldots,a_r^m)_v
\subseteq s^q d^mD\subseteq(a_1^m,\ldots,a_r^m)_v.
\qedhere
\]
\end{proof}

The power estimate gives the finite-family form of the definition and
allows all the exponents occurring in a finite collection of witnesses
to be synchronized.

\begin{theorem}\label{t:finite-families}
Let $D$ be a domain and let $S\subseteq D^\ast$ be multiplicative.
The following conditions are equivalent.
\begin{enumerate}
\item $D$ is almost $S$-GCD.
\item For every finite nonempty family $a_1,\ldots,a_r\in D^\ast$,
there exists $n\geq1$ such that $(a_1^n,\ldots,a_r^n)_v$ is
$S$-principal.
\item For every finite nonempty family $a_1,\ldots,a_r\in D^\ast$,
there exists $n\geq1$ such that $\bigcap_{i=1}^r a_i^nD$ is
$S$-principal.
\end{enumerate}
Moreover, an exponent witnessing either of the last two conditions
can be replaced by any positive multiple.
\end{theorem}

\begin{proof}
We prove the implication from the first condition to the second by
induction on $r$. The case $r=1$ is immediate. Suppose that
$s(a_1^n,\ldots,a_r^n)_v\subseteq dD\subseteq(a_1^n,\ldots,a_r^n)_v$
for some $s\in S$ and $d\in D^\ast$.
Apply the defining property to $d$ and $a_{r+1}^n$ to choose $m\geq1$,
$t\in S$, and $e\in D^\ast$ with
$tK\subseteq eD\subseteq K$, where
$K=(d^m,a_{r+1}^{nm})_v$.
Put $J=(a_1^{nm},\ldots,a_{r+1}^{nm})_v$ and
$q=(r-1)(m-1)$. Lemma~\ref{l:power-sandwich} gives
$s^mJ\subseteq K$ and $s^qK\subseteq J$. Hence
$s^{m+q}tJ\subseteq s^qeD\subseteq J$.
This proves the induction step.

Inversion and multiplication by a nonzero scalar preserve
$S$-principality for divisorial fractional ideals, as shown in the
proof of Theorem~\ref{t1}. If $A=\prod_i a_i$, then
\[
\bigcap_i a_i^nD
=A^n\bigl((A/a_1)^n,\ldots,(A/a_r)^n\bigr)^{-1}.
\]
This proves that the second condition implies the third. Applying
the third condition to $A/a_1,\ldots,A/a_r$ proves the reverse
implication, since
$\bigcap_i(A/a_i)^nD=A^n(a_1^n,\ldots,a_r^n)^{-1}$.
Taking two elements and using Theorem~\ref{t1} gives the first
condition. The assertion about multiples follows from
Lemma~\ref{l:power-sandwich} and the same inversion identities.
\end{proof}

The following permanence properties clarify the dependence on the
multiplicative set. Its saturation is
$\Sat(S)=\{x\in D^\ast:xy\in S\text{ for some }y\in D\}$.

\begin{proposition}\label{p:permanence}
\label{p:enlarging-S}\label{p:localization-stability}\label{t2}
Let $D$ be a domain and let $S,T\subseteq D^\ast$ be multiplicative.
\begin{enumerate}
\item If $S\subseteq T$ and $D$ is almost $S$-GCD, then $D$ is
almost $T$-GCD.
\item $D$ is almost $S$-GCD if and only if it is almost
$\Sat(S)$-GCD.
\item If $D$ is almost $S$-GCD, then $D_T$ is almost
$S_T$-GCD, where $S_T=\{s/1:s\in S\}$.
In particular, $D_S$ is AGCD.
\end{enumerate}
Consequently, an AGCD domain is almost $S$-GCD for every $S$.
\end{proposition}

\begin{proof}
The first assertion follows by retaining an $S$-principal witness.
For the second, only descent from $\Sat(S)$ needs proof.
If $uI\subseteq dD\subseteq I$ with $u\in\Sat(S)$, choose
$v\in D$ with $vu\in S$. Then
$(vu)I\subseteq vdD\subseteq dD\subseteq I$, so $I$ is
$S$-principal.

For the third, write two nonzero elements of $D_T$ as $a/u$ and
$b/v$, with $a,b\in D^\ast$ and $u,v\in T$. Choose
$n\geq1$ such that $a^nD\cap b^nD$ is $S$-principal and localize
its defining sandwich. Localization commutes with finite
intersections, and the units $u,v$ do not affect the corresponding
principal ideals. Thus
$(a/u)^nD_T\cap(b/v)^nD_T$ is $S_T$-principal.
For $T=S$, every element of $S_T$ is a unit, giving the AGCD
condition. The final assertion follows from the first with the
smaller set $\{1\}$.
\end{proof}

For a $D$-module $E$, write
$\tau_S(E)=\ker(E\longrightarrow E_S)
=\{x\in E:sx=0\text{ for some }s\in S\}$.
The module $E$ is \emph{$S$-torsion} if $\tau_S(E)=E$, and it is
\emph{uniformly $S$-torsion} if $sE=0$ for one $s\in S$.
Uniform torsion records the single denominator required by
$S$-principality.

\begin{proposition}\label{p:torsion-obstruction}
Let $D$ be a domain and let $S\subseteq D^\ast$ be multiplicative.
\begin{enumerate}
\item A nonzero fractional ideal $I$ is $S$-principal if and only
if $I/dD$ is uniformly $S$-torsion for some $d\in I\setminus\{0\}$.
If $ID_S$ is principal and $d\in I$ satisfies $ID_S=dD_S$, then
$I/dD$ is $S$-torsion. In this case, $I$ is $S$-principal if and
only if $I/dD$ is uniformly $S$-torsion; this criterion is
independent of the chosen generator $d$ of $ID_S$.
\item For each $d\in D^\ast$, there is a natural identification
\[
\tau_S(D/dD)=(dD_S\cap D)/dD.
\]
Thus the condition
\begin{equation}\tag{$C_{D,S}$}\label{eq:C}
\text{for every }d\in D^\ast\text{ there exists }s\in S
\text{ such that }s(dD_S\cap D)\subseteq dD
\end{equation}
holds if and only if $\tau_S(D/dD)$ is uniformly $S$-torsion
for every $d\in D^\ast$. Equivalently, for every $d\in D^\ast$
there is $s\in S$ with $dD_S\cap D=(dD:s)$.
\item If $S=\langle u\rangle$ for $u\in D^\ast$, then
\eqref{eq:C} holds if and only if, for every $d\in D^\ast$,
the ascending chain $((dD:u^j))_{j\geq0}$ eventually stabilizes.
\end{enumerate}
\end{proposition}

\begin{proof}
The first equivalence in (1) is the definition of the quotient
annihilator. If $ID_S=dD_S$, then each element of $I/dD$ is killed
by some element of $S$. Suppose that $I$ is $S$-principal, with
$sI\subseteq eD\subseteq I$. Since $e\in dD_S$, choose $t\in S$
such that $te\in dD$. Then $tsI\subseteq dD$, proving uniform
torsion for the specified $d$ and hence independence of that choice.

An element $x+dD$ belongs to $\tau_S(D/dD)$ precisely when
$tx\in dD$ for some $t\in S$, which proves the identification in
(2). If $s(dD_S\cap D)\subseteq dD$, then
$dD_S\cap D\subseteq(dD:s)$; the reverse inclusion always holds.
This proves all the equivalences in (2). Finally,
$dD_u\cap D=\bigcup_{j\geq0}(dD:u^j)$, so the last assertion
follows from (2).
\end{proof}

The difference between elementwise and uniform torsion also has a
spectral description. For an ideal $J$, put
$V(J)=\{P\in\Spec(D):J\subseteq P\}$, and set
$U_S=\{P\in\Spec(D):P\cap S=\varnothing\}$.
This subset is the image of $\Spec(D_S)$ in $\Spec(D)$; it need
not be Zariski open. For a module $E$, its support is
$\Supp_D(E)=\{P\in\Spec(D):E_P\neq0\}$.

\begin{proposition}\label{p:annihilator-topology}
Let $D$ be a domain, let $S\subseteq D^\ast$ be multiplicative,
and let $E$ be a $D$-module. Then:
\begin{enumerate}
\item $E$ is $S$-torsion if and only if
$\Supp_D(E)\cap U_S=\varnothing$;
\item $E$ is uniformly $S$-torsion if and only if
$V(\Ann_D(E))\cap U_S=\varnothing$.
\end{enumerate}
If $E$ is finitely generated, these two conditions are equivalent.
\end{proposition}

\begin{proof}
Localization identifies $\Supp_{D_S}(E_S)$ with
$\Supp_D(E)\cap U_S$. A module is zero if and only if all of
its localizations at prime ideals are zero. This proves (1).
For (2), the stated intersection is empty precisely when
$\Ann_D(E)D_S=D_S$, or equivalently when
$\Ann_D(E)\cap S\neq\varnothing$.
If $E$ is finitely generated and $E_S=0$, choose an element of $S$
annihilating each member of a finite generating set and multiply
these elements. Their product annihilates $E$.
\end{proof}

The finiteness qualification in the last proposition is necessary.

\begin{example}\label{e:unbounded-torsion}
Let $D=\ZZ$, let $p$ be a prime number, put
$S=\langle p\rangle$, and take $E=\ZZ[1/p]/\ZZ$.
Every element of $E$ is killed by a power of $p$, so $E_S=0$.
However, $p^j(1/p^{j+1}+\ZZ)\neq0$ for every $j\geq0$;
thus $E$ is not uniformly $S$-torsion. Moreover,
$\Ann_{\ZZ}(E)=0$ and $\Supp_{\ZZ}(E)=\{(p)\}$.
Consequently, $\Supp_{\ZZ}(E)$ avoids $U_S$, whereas
$V(\Ann_{\ZZ}(E))=\Spec(\ZZ)$ does not.
\end{example}

The following descent theorem is the power analogue of
\cite[Theorem~2.3]{AHZ}. The module formulation above specifies
exactly what its contraction hypothesis requires.

\begin{theorem}\label{t3}
Let $D$ be a domain and let $S\subseteq D^\ast$ be multiplicative.
Assume that $(C_{D,S})$ holds. Then $D$ is almost $S$-GCD if and
only if $D_S$ is AGCD.
\end{theorem}

\begin{proof}
The forward implication follows from Proposition~\ref{p:permanence}.
Conversely, take $a,b\in D^\ast$. Since $D_S$ is AGCD, choose
$n\geq1$ such that $(a^nD\cap b^nD)D_S$ is principal.
Its generator can be chosen in $a^nD\cap b^nD$ after multiplication
by a unit from $S$; write it as $d$. Then
$a^nD\cap b^nD\subseteq dD_S\cap D$.
By $(C_{D,S})$, there exists $s\in S$ with
\[
s(a^nD\cap b^nD)\subseteq s(dD_S\cap D)
\subseteq dD\subseteq a^nD\cap b^nD.
\]
Theorem~\ref{t1} now gives the result.
\end{proof}

A finiteness condition on the particular ideals being descended can
replace the stronger requirement on all principal contractions.

\begin{lemma}\label{l:finite-obstruction}
Let $D$ be a domain, let $S\subseteq D^\ast$ be multiplicative,
and let $I$ be a nonzero fractional divisorial ideal of $D$.
Then $I$ is $S$-principal if and only if $ID_S$ is principal and
there exist $u\in S$ and a finitely generated fractional ideal
$F\subseteq I$ such that $uI\subseteq F_v\subseteq I$.
\end{lemma}

\begin{proof}
The forward implication follows by taking $F=dD$ in an
$S$-principal sandwich. Conversely, choose $d\in I$ with
$ID_S=dD_S$. Since $F$ is finitely generated, some $v\in S$
satisfies $vF\subseteq dD$. Taking $v$-closures yields
$vF_v\subseteq dD$, so $uvI\subseteq dD\subseteq I$.
\end{proof}

\begin{theorem}\label{t:finite-conductor-localization}
Let $D$ be a domain and let $S\subseteq D^\ast$ be multiplicative.
Suppose that, for every $a,b\in D^\ast$, the ideal $J=aD\cap bD$
satisfies $uJ\subseteq F_v\subseteq J$ for some $u\in S$ and
some finitely generated ideal $F\subseteq J$. Then $D$ is almost
$S$-GCD if and only if $D_S$ is AGCD.
\end{theorem}

\begin{proof}
Only descent needs proof. For $a,b\in D^\ast$, choose $n\geq1$
such that $(a^nD\cap b^nD)D_S$ is principal. The divisorial ideal
$a^nD\cap b^nD$ satisfies the stated finiteness hypothesis, so
Lemma~\ref{l:finite-obstruction} applies.
\end{proof}

Recall that a ring is coherent if every finitely generated ideal is
finitely presented. The preceding theorem applies in particular to
the usual finiteness hypotheses on intersections of principal ideals.

\begin{corollary}\label{c:finite-conductor}
Let $D$ be a domain and let $S\subseteq D^\ast$ be multiplicative.
If every $aD\cap bD$, with $a,b\in D^\ast$, is of $v$-finite type,
then $D$ is almost $S$-GCD if and only if $D_S$ is AGCD.
In particular, this equivalence holds when every such intersection
is finitely generated, and hence when $D$ is coherent.
\end{corollary}

\begin{proof}
Use $u=1$ in Theorem~\ref{t:finite-conductor-localization}.
For the coherent case, the kernel of the surjection
$D^2\longrightarrow aD+bD$, $(x,y)\mapsto ax+by$, is finitely
generated and is isomorphic to $aD\cap bD$ via $(x,y)\mapsto ax$.
\end{proof}

A global finite-type hypothesis gives a direct contraction proof.
Recall that a Mori domain satisfies the ascending chain condition
on integral divisorial ideals.

\begin{corollary}\label{c:tfinite-localization-criterion}
Let $D$ be a domain and let $S\subseteq D^\ast$ be multiplicative.
If every nonzero integral $t$-ideal of $D$ is of $v$-finite type,
then $(C_{D,S})$ holds and $D$ is almost $S$-GCD if and only if
$D_S$ is AGCD. This applies in particular to Mori domains, and
hence to Noetherian domains.
\end{corollary}

\begin{proof}
For $d\in D^\ast$, put
$K=dD_S\cap D=\bigcup_{s\in S}(dD:s)$.
Each $(dD:s)=s^{-1}dD\cap D$ is divisorial. The union is directed,
so the finite-type property of $t$ shows that $K$ is a $t$-ideal.
By hypothesis, $K=F_v$ for a finitely generated $F\subseteq K$.
Choose $s\in S$ with $sF\subseteq dD$. Then
$sK=(sF)_v\subseteq dD$, proving $(C_{D,S})$.
Apply Theorem~\ref{t3}.

If $D$ is Mori, every nonzero integral ideal $L$ has a finitely
generated subideal $F$ with $F_v=L_v$: otherwise one constructs
a strictly ascending chain of the $v$-closures of finite subideals
of $L$. If $L$ is a $t$-ideal, then $F_v\subseteq L$, forcing
$L=F_v$. Thus Mori domains satisfy the stated hypothesis.
\end{proof}

A saturated multiplicative set $S\subseteq D^\ast$ is a
\emph{splitting set} if every $d\in D^\ast$ can be written
$d=sa$ with $s\in S$, $a\in D^\ast$, and
$uD\cap aD=uaD$ for every $u\in S$ \cite{AAZ}.

\begin{corollary}\label{c:splitting-localization-criterion}
Let $D$ be a domain and let $S\subseteq D^\ast$ be a splitting
set. Then $(C_{D,S})$ holds, and $D$ is almost $S$-GCD if and
only if $D_S$ is AGCD.
\end{corollary}

\begin{proof}
Write $d=sa$ as in the definition. We claim that
$dD_S\cap D=aD$. Since $s$ is a unit in $D_S$, it suffices to
show $aD_S\cap D=aD$. If $x\in aD_S\cap D$, choose $u\in S$
with $ux\in aD$. Then $ux\in uD\cap aD=uaD$, so $x\in aD$.
The reverse inclusion is immediate. Hence
$s(dD_S\cap D)=saD=dD$, and Theorem~\ref{t3} applies.
\end{proof}

The following example strengthens the localization counterexample
from \cite[Example~2.2]{AHZ}: all positive powers of the indicated
pair fail, and the domain is integrally closed.

\begin{example}\label{e2}
Let $p$ be a prime number, put $A=\QQ[[Y]]$ and
$V=\ZZ_{(p)}+YA$, and set
$D=V+XA[X]$ and $S=\langle p\rangle$.
Then $D$ is integrally closed and $D_S=A[X]$ is a UFD, but
$D$ is not almost $S$-GCD.

Indeed, $V$ is the valuation domain of $\QQ((Y))$ associated to
the lexicographically ordered valuation
$z\mapsto(\operatorname{ord}_Y(z),v_p(c_z))$, where $c_z$ is
the first nonzero coefficient of $z$. Thus $V$ is integrally
closed. Also $V_S=A$, whence $D_S=A[X]$.
If $z\in\qf(D)$ is integral over $D$, it is integral over the
integrally closed domain $A[X]$ and therefore lies in $A[X]$.
Evaluating a monic equation for $z$ over $D$ at $X=0$ shows
that $z(0)$ is integral over $V$. Thus $z(0)\in V$ and
$z\in V+XA[X]=D$, proving normality.

For every $n\geq1$, one has
\[
J_n=X^nD\cap Y^nD=X^nY^nA[X].
\]
To see the nontrivial inclusion, write $f=X^ng=Y^nh$ with
$g,h\in D$. In $A[X]$, divisibility by $X^n$ implies
$h=X^nk$ for some $k\in A[X]$, and hence $f=X^nY^nk$.
Conversely, both $Y^nk$ and $X^nk$ belong to $D$ for every
$k\in A[X]$. Since $p$ is a unit in $A$, we have $pJ_n=J_n$.
If $J_n$ were $S$-principal, a sandwich
$p^jJ_n\subseteq dD\subseteq J_n$ would give $J_n=dD$.
It would then follow that $pdD=dD$, so $p$ would be a unit
of $D$. This is impossible because evaluation at $X=0$ would
give $1/p\in V$. Theorem~\ref{t1} proves the assertion.
\end{example}

The principal-contraction condition is sufficient for descent but
is not necessary, even in a valuation domain.

\begin{example}\label{e:C-not-necessary-valuation}
Let $p$ be a prime number, let
$D=\ZZ_{(p)}+Y\QQ[[Y]]$, and let $S=\langle p\rangle$.
The domain $D$ is a valuation domain, so it is GCD and hence
almost $S$-GCD. Also $D_S=\QQ[[Y]]$ is a DVR. Nevertheless,
$(C_{D,S})$ fails. In fact,
\[
YD_S\cap D=Y\QQ[[Y]],\qquad
YD=Y\ZZ_{(p)}+Y^2\QQ[[Y]].
\]
The first ideal is unchanged by multiplication by any element
of $S$ and strictly contains the second: $Y/p$ belongs to the
first but not the second. Thus no $s\in S$ satisfies
$s(YD_S\cap D)\subseteq YD$.
\end{example}

\section{Relative class groups and divisor classes}\label{sec:classes}

The relative class group records the obstruction to making a
$t$-invertible ideal principal after multiplying by a denominator in
$S$. We first identify its localization map and its kernel. The
resulting torsion criterion applies to almost Pr\"ufer
$v$-multiplication domains without a bounded-contraction hypothesis.
For Krull domains, divisor theory makes this criterion explicit in
terms of the height-one primes meeting $S$.

\subsection{Localization and torsion}

Let $T(D)$ be the group of $t$-invertible fractional $t$-ideals of
$D$, with product $I\cdot_t J=(IJ)_t$, and let $P(D)$ be its subgroup
of nonzero principal fractional ideals. Thus
$\Cl_t(D)=T(D)/P(D)$. Following \cite{HamedNagata, Hamed}, put
\[
 S\text{-}P(D)=\{I\in T(D):I\text{ is }S\text{-principal}\},
 \qquad
 S\text{-}\Cl_t(D)=T(D)/(S\text{-}P(D)).
\]
The first set is a subgroup: if $sI\subseteq aD\subseteq I$ and
$uJ\subseteq bD\subseteq J$, taking $t$-closures of their products
gives $su(IJ)_t\subseteq abD\subseteq(IJ)_t$; inversion gives
$sI^{-1}\subseteq(s/a)D\subseteq I^{-1}$. It contains $P(D)$,
so the relative group is also a quotient of $\Cl_t(D)$.

The next lemma supplies a proof of the localization facts needed
from this relative class-group framework.

\begin{lemma}\label{l:class-localization-kernel}
Let $D$ be a domain, let $S\subseteq D^\ast$ be multiplicative,
and let $I\in T(D)$. Then $ID_S$ is a $t$-invertible divisorial
fractional ideal of $D_S$. It is principal if and only if $I$ is
$S$-principal. Localization therefore induces an injective homomorphism
\[
 S\text{-}\Cl_t(D)\longrightarrow\Cl_t(D_S).
\]
If $D$ is a PvMD, this homomorphism is an isomorphism.
\end{lemma}

\begin{proof}
Choose finitely generated fractional ideals $F,G$ with
$I=F_v$ and $I^{-1}=G_v$. Inversion of a finitely generated ideal
is a finite intersection of principal fractional ideals and commutes
with localization. Consequently,
\[
 ID_S=(G^{-1})D_S=(GD_S)^{-1},\qquad
 I^{-1}D_S=(FD_S)^{-1}.
\]
Here and throughout this proof, inverses of $D_S$-ideals are taken
with respect to $D_S$. Put $J=I^{-1}D_S$. Since
$(ID_S)J\subseteq D_S$ and $GD_S\subseteq J$, we have
\[
 ID_S\subseteq J^{-1}\subseteq(GD_S)^{-1}=ID_S.
\]
Thus $ID_S=J^{-1}=(FD_S)_{v_{D_S}}$ is divisorial, and its inverse
is $J$. Choose a finitely generated ideal $H\subseteq II^{-1}$
with $H_v=D$. Since $H^{-1}=D$, localization gives
$(HD_S)^{-1}=D_S$, whence $(HD_S)_{v_{D_S}}=D_S$.
It follows that $((ID_S)J)_{t_{D_S}}=D_S$.

For $I=F_v,L=E_v\in T(D)$ with $F,E$ finitely generated,
their $t$-product is $(FE)_v$ and belongs to $T(D)$. Applying
the identity just proved to this product gives
\[
 ((IL)_t)D_S=(FED_S)_{v_{D_S}}
            =((ID_S)(LD_S))_{t_{D_S}}.
\]
Localization is therefore a homomorphism on the groups of
$t$-invertible $t$-ideals.

If $sI\subseteq aD\subseteq I$, localization gives $ID_S=aD_S$.
Conversely, if $ID_S=cD_S$, multiply $c$ by an element of $S$
to obtain $a\in I$ with $ID_S=aD_S$. Since $F$ is finitely
generated, some $s\in S$ satisfies $sF\subseteq aD$.
Taking $v$-closures gives $sI\subseteq aD\subseteq I$.
Thus the kernel of $T(D)\to\Cl_t(D_S)$ is exactly $S\text{-}P(D)$,
which proves injectivity of the stated map.

Finally suppose $D$ is a PvMD, and let $L\in T(D_S)$.
Write $L=(a_1D_S+\cdots+a_rD_S)_{v_{D_S}}$ with $a_i\in K^\ast$.
Then $I=(a_1D+\cdots+a_rD)_v$ is $t$-invertible, and the
preceding localization identity gives $ID_S=L$. The map is
therefore surjective.
\end{proof}

The following is the finite fractional-ideal form of the power
lemma of Anderson and Zafrullah \cite[Lemma~3.3]{AZ-Bez}.
We recall it for use in the relative class-group argument.

\begin{lemma}\label{l:t-invertible-powers}
Let $D$ be a domain with quotient field $K$, and let
$a_1,\ldots,a_r\in K\setminus\{0\}$. Assume that
$F=(a_1,\ldots,a_r)D$ is $t$-invertible. Then, for every
integer $n\geq1$,
\[
 ((F_t)^n)_t=(F^n)_t
 =(a_1^n,\ldots,a_r^n)_t
 =(a_1^n,\ldots,a_r^n)_v.
\]
Moreover, the common fractional $t$-ideal is $t$-invertible.
\end{lemma}

\begin{proof}
Choose $0\ne h\in D$ such that $ha_i\in D$ for every $i$.
Since $(hF)_t=hF_t$ is $t$-invertible,
\cite[Lemma~3.3]{AZ-Bez} gives
\[
\begin{aligned}
 h^n(F^n)_t
 &=((hF)^n)_t \\
 &=((ha_1)^n,\ldots,(ha_r)^n)_t \\
 &=h^n(a_1^n,\ldots,a_r^n)_t.
\end{aligned}
\]
Cancelling $h^n$ yields
$(F^n)_t=(a_1^n,\ldots,a_r^n)_t$.
The identity $((F_t)^n)_t=(F^n)_t$ follows from the
multiplicative properties of the $t$-operation.

The ideal generated by $a_1^n,\ldots,a_r^n$ is finitely generated,
so its $t$-closure equals its $v$-closure.
Finally, $(F^n)_t$ is $t$-invertible because $F$ is
$t$-invertible.
\end{proof}

The torsion conclusion holds already under the localized AGCD
condition, whether or not descent of almost $S$-GCD is available.

\begin{theorem}\label{S-class group is torsion}
\label{c:relative-torsion-weaker-hypothesis}
Let $D$ be a domain and let $S\subseteq D^\ast$ be multiplicative.
If $D_S$ is AGCD, then $S\text{-}\Cl_t(D)$ is torsion.
In particular, this conclusion holds if $D$ is almost $S$-GCD.
\end{theorem}

\begin{proof}
The $t$-class group of an AGCD domain is torsion
\cite[Theorem~3.4]{AZ-Bez}. Apply the injection of
Lemma~\ref{l:class-localization-kernel}.
For the final assertion, Theorem~\ref{t2} shows that $D_S$ is AGCD.
\end{proof}

The APvMD hypothesis ensures that a suitable power of each pair
defines a class to which the torsion condition can be applied.
We use the notion introduced and studied in \cite{LiAPvMD}.

\begin{samepage}
\begin{theorem}\label{t:apvmd-relative-characterization}
\label{t:pvmd-relative-characterization}
Let $D$ be an APvMD and let $S\subseteq D^\ast$ be multiplicative.
The following conditions are equivalent.
\begin{enumerate}
\item $D$ is almost $S$-GCD.
\item $D_S$ is AGCD.
\item $S\text{-}\Cl_t(D)$ is torsion.
\end{enumerate}
In particular, the equivalence holds for every PvMD, without a
principal-contraction hypothesis.
\end{theorem}
\end{samepage}

\begin{proof}
The first condition implies the second by Theorem~\ref{t2},
and the second implies the third by
Theorem~\ref{S-class group is torsion}.
Assume the third condition and take $a,b\in D^\ast$.
Since $D$ is an APvMD, there is $m\geq1$ such that
$F=(a^m,b^m)D$ is $t$-invertible. Put $I=F_t\in T(D)$.
The class of $I$ in $S\text{-}\Cl_t(D)$ has finite order,
so $(I^n)_t$ is $S$-principal for some $n\geq1$.
Lemma~\ref{l:t-invertible-powers} gives
$(I^n)_t=(a^{mn},b^{mn})_v$. Thus $D$ is almost $S$-GCD.
\end{proof}

The hypotheses and conclusions are separated by the examples from
the preceding section.

\begin{remark}\label{r:class-torsion-limits}
In Example~\ref{e2}, the domain $D$ is integrally closed and $D_S$
is a UFD, but $D$ is not almost $S$-GCD. The injection in
Lemma~\ref{l:class-localization-kernel} gives
$S\text{-}\Cl_t(D)=0$. Thus even triviality of the relative class
group is insufficient for an arbitrary integrally closed domain;
the APvMD assumption cannot be replaced merely by normality.
Conversely, Example~\ref{e:C-not-necessary-valuation} shows that
the bounded-contraction condition can fail in a GCD domain.
The APvMD characterization therefore applies beyond that
descent condition.
\end{remark}

\subsection{Divisor classes in Krull domains}

For a Krull domain $D$, write $X^1(D)$ for its height-one prime
ideals and $\Cl(D)$ for its divisor class group. In this setting
all nonzero fractional divisorial ideals are $t$-invertible and
$\Cl_t(D)=\Cl(D)$. The following formulation combines the
preceding characterization with the classical localization
sequence for divisor classes; see \cite{Fossum}.

\begin{theorem}\label{t:krull-boundary-criterion}
Let $D$ be a Krull domain, let $S\subseteq D^\ast$ be
multiplicative, and put
\[
 H_S=\langle[P]:P\in X^1(D),\ P\cap S\ne\varnothing\rangle
       \subseteq\Cl(D).
\]
Then
\[
 S\text{-}\Cl_t(D)\simeq\Cl(D_S)\simeq\Cl(D)/H_S.
\]
Consequently, $D$ is almost $S$-GCD if and only if
$\Cl(D)/H_S$ is torsion. Equivalently, the vectors $[P]\otimes1$
for the height-one primes meeting $S$ span
$\Cl(D)\otimes_{\ZZ}\QQ$.
\end{theorem}

\begin{proof}
A Krull domain is a PvMD, so the first isomorphism follows from
Lemma~\ref{l:class-localization-kernel}. Its divisor group is the
free abelian group on $X^1(D)$. Localization deletes the basis
elements indexed by primes meeting $S$ and retains those indexed
by primes disjoint from $S$. It sends the divisor of an element
$x\in K^\ast$ to the divisor of the same element in $D_S$.
The resulting class-group map is surjective, since any finite
divisor of $D_S$ lifts to a finite divisor of $D$.
If the image of a divisor $A$ is the divisor of $x$ in $D_S$,
then $A-\operatorname{div}_D(x)$ is supported on the height-one
primes meeting $S$. This identifies the kernel as $H_S$ and
proves the second isomorphism.

The almost $S$-GCD criterion now follows from
Theorem~\ref{t:apvmd-relative-characterization}.
An abelian group is torsion precisely when its tensor product
with $\QQ$ is zero. Tensoring the quotient by $\QQ$ yields the
stated spanning criterion.
\end{proof}

Finite rational class rank permits all needed prime-divisor classes
to be removed by one denominator.

\begin{corollary}\label{c:one-denominator-finite-rank}
Let $D$ be a Krull domain such that
$\dim_{\QQ}(\Cl(D)\otimes_{\ZZ}\QQ)<\infty$, and let
$S\subseteq D^\ast$ be multiplicative. Then $D$ is almost
$S$-GCD if and only if it is almost $\langle s\rangle$-GCD
for some $s\in S$.
\end{corollary}

\begin{proof}
One implication follows from $\langle s\rangle\subseteq S$.
For the converse, Theorem~\ref{t:krull-boundary-criterion}
provides finitely many primes $P_1,\ldots,P_r$ meeting $S$
whose classes span $\Cl(D)\otimes_{\ZZ}\QQ$.
Choose $s_i\in S\cap P_i$ and put $s=s_1\cdots s_r$.
Every $P_i$ contains $s$, so the same theorem applied to
$\langle s\rangle$ proves the conclusion.
If the rational class group is zero, take $s=1$.
\end{proof}

A bound on the exponent of the relative class group yields a
common exponent for all finite families.

\begin{corollary}\label{c:uniform-class-exponent}
Let $D$ be a PvMD, let $S\subseteq D^\ast$ be multiplicative,
and let $e\geq1$. The following conditions are equivalent.
\begin{enumerate}
\item $e(S\text{-}\Cl_t(D))=0$.
\item For every finite nonempty family $a_1,\ldots,a_r\in K^\ast$,
the ideal $(a_1^e,\ldots,a_r^e)_v$ is $S$-principal.
\end{enumerate}
In particular, if $D$ is Krull, $\Cl(D)/H_S$ is finitely
generated, and $D$ is almost $S$-GCD, then one exponent $e$
works simultaneously for all such finite families.
\end{corollary}

\begin{proof}
In a PvMD, $I=(a_1,\ldots,a_r)_v$ belongs to $T(D)$, and
Lemma~\ref{l:t-invertible-powers} identifies $(I^e)_t$ with
$(a_1^e,\ldots,a_r^e)_v$. Conversely, every member of $T(D)$
can be written in this form. These observations
prove the equivalence. In the final assertion,
Theorem~\ref{t:krull-boundary-criterion} makes
$\Cl(D)/H_S$ a finitely generated torsion abelian group.
It is finite, so its exponent supplies the required $e$.
\end{proof}

The distinction between almost $S$-GCD and $S$-GCD already
occurs in a Noetherian Krull domain with finite class group.

\begin{example}\label{e1}
Let $p$ be a rational prime and $n\geq2$. Put
\[
 D=\ZZ[X,Y]/(XY-p^n)\simeq\ZZ[X,p^n/X],
 \qquad M=(p,X,Y)D,
 \qquad S=D\setminus M,
\]
where the residue classes of $X,Y$ retain their names.
The domain $D$ is Krull and $\Cl(D)\simeq\ZZ/n\ZZ$, as in
\cite{AA2}. These facts have a direct verification here.
Inside $\QQ(X)$, with $Y=p^n/X$, one has
\[
 D=D_X\cap D_Y
   =\ZZ[X,X^{-1}]\cap\ZZ[Y,Y^{-1}].
\]
Indeed, a Laurent polynomial $\sum_i c_iX^i$ in the first ring
belongs to the second precisely when $c_i p^{ni}\in\ZZ$ for
each $i$. For $i<0$, this means that $c_i$ is divisible by
$p^{-ni}$, which is exactly the description of $D$.
The intersection of these two integrally closed rings is
integrally closed; since $D$ is Noetherian, it is Krull.

The only height-one prime containing $X$ is $P=(p,X)D$.
In $D_P$, the element $Y$ is a unit and $X=p^n/Y$, so $p$
is a uniformizer and $\operatorname{div}_D(X)=nP$.
Since $D_X$ is a UFD, the divisor localization sequence in
Theorem~\ref{t:krull-boundary-criterion} makes $\Cl(D)$
generated by $[P]$. If $kP=\operatorname{div}_D(f)$, then $f$
is a unit of $D_X$, hence $f=\pm X^j$ for some $j\in\ZZ$.
Thus $k=nj$, proving that $[P]$ has exact order $n$.
Theorem~\ref{t:apvmd-relative-characterization}, with
multiplicative set $\{1\}$, now shows that $D$ is AGCD and
consequently almost $S$-GCD.

However, $D_M$ is not a UFD. To see this directly, let
$\mathfrak m=MD_M$. The element $p$ lies in
$\mathfrak m\setminus\mathfrak m^2$: evaluation at $X=Y=0$
defines a homomorphism
$D_M\to\ZZ_{(p)}/p^n\ZZ_{(p)}$, under which the image of
$\mathfrak m^2$ is contained in the square of $(p)$, whereas
the image of $p$ is not in that square. Hence $p$ cannot
factor as a product of two nonunits of $D_M$ and is irreducible.
It is not prime, since
\[
 D_M/pD_M\simeq
 \bigl(\mathbb F_p[X,Y]/(XY)\bigr)_{(X,Y)}
\]
has nonzero zero divisors, namely the images of $X$ and $Y$.
Since a Noetherian GCD domain is a UFD, $D_M=D_S$ is not GCD.
Localization of an $S$-GCD domain at $S$ is GCD
\cite[Theorem~2.2]{AHZ}. Therefore $D$ is not $S$-GCD.
\end{example}

\section{Normalization and polynomial extensions}\label{sec:polynomial}

Localization exposes the classical root-extension structure of AGCD
domains. We first use it to describe normalization and its spectrum away
from $S$. We then study polynomial extensions, distinguishing the
additional root-extension obstruction from the bounded-denominator
condition. Over an integrally closed domain, the latter condition is
already determined by the coefficient domain. Almost splitting sets
provide another setting in which polynomial transfer is possible.

\subsection{Root extensions and normalization}

We begin with the standard normalization theorem for AGCD domains,
applied after localization.

\begin{proposition}\label{p:localized-root-extension}
Let $D$ be a domain and let $S\subseteq D^\ast$ be multiplicative.
If $D_S$ is an AGCD domain, then
$D_S\subseteq\overline{D_S}=(\overline D)_S$ is a root extension,
and $(\overline D)_S$ is an integrally closed AGCD domain.
In particular, for every $x\in\overline D$ there are $n\geq1$ and
$s\in S$ such that $sx^n\in D$. These conclusions apply whenever
$D$ is almost $S$-GCD.
\end{proposition}
\begin{proof}
The integral closure of an AGCD domain is again AGCD, and the
normalization extension is a root extension
\cite[Theorems~3.1 and~3.4]{Z3}. Integral closure commutes with
localization. Thus each $x\in\overline D$ has a power in $D_S$,
and clearing its denominator gives the assertion.
The final statement follows from Theorem~\ref{t2}.
\end{proof}

For $x\in K$, call $x$ \emph{$S$-integral over $D$} if $sx$ is
integral over $D$ for some $s\in S$, and write $\overline D^{\,S}$
for the set of such elements. This ring is always $(\overline D)_S$;
under the AGCD hypothesis it also has the following power description.

\begin{theorem}\label{t:S-integral-root-test}
Let $D$ be a domain, let $S\subseteq D^\ast$ be multiplicative,
and assume that $D_S$ is an AGCD domain. Then
\[
 \overline D^{\,S}=(\overline D)_S=\overline{D_S}
 =\{x\in K: sx^m\in D\text{ for some }s\in S, m\geq1\}.
\]
\end{theorem}
\begin{proof}
The first equality follows directly from the definition, and the
second is compatibility of normalization with localization.
If $x\in\overline{D_S}$, Proposition~\ref{p:localized-root-extension}
gives $x^m\in D_S$ for some $m\geq1$, so $sx^m\in D$ for some
$s\in S$. Conversely, $sx^m\in D$ implies
$(sx)^m=s^{m-1}(sx^m)\in D$. Hence $sx$ is integral over $D$,
so $x\in\overline D^{\,S}$.
\end{proof}

Recall that $D$ is \emph{root closed} if $x\in K$ and $x^n\in D$
for some $n\geq1$ imply $x\in D$.

\begin{corollary}\label{c:root-closed-localization}
Let $D$ be a domain and $S\subseteq D^\ast$ multiplicative.
If $D_S$ is AGCD and root closed, then $D_S$ is integrally closed.
In particular, this holds if $D$ is root closed and almost $S$-GCD.
\end{corollary}
\begin{proof}
A root extension from a root-closed domain to its integral closure
is the identity. For the final assertion, root closedness passes
to localization: if $x^n=a/s\in D_S$, then
$(sx)^n=s^{n-1}a\in D$, so $sx\in D$ and $x\in D_S$.
\end{proof}

\begin{corollary}\label{c:integral-elements-cleared}
Let $D$ be a domain and $S\subseteq D^\ast$ multiplicative.
If $D_S$ is AGCD and root closed, then every $x\in\overline D$
satisfies $sx\in D$ for some $s\in S$.
\end{corollary}
\begin{proof}
The preceding corollary gives
$\overline D\subseteq\overline{D_S}=D_S$.
\end{proof}

For a domain $A$ containing $D$, put
$U_S(A)=\{P\in\Spec(A):P\cap S=\varnothing\}$, with the
subspace topology. For arbitrary $S$ this subspace need not be open.
The root-extension property determines the primes of the normalization
lying over this subspace.

\begin{theorem}\label{t:normalization-homeomorphism}
Let $D$ be a domain, let $S\subseteq D^\ast$ be multiplicative,
and assume that $D_S$ is AGCD. Contraction induces a homeomorphism
$U_S(\overline D)\longrightarrow U_S(D)$.
For $P\in U_S(D)$, its unique inverse image is
\[
 Q(P)=\{x\in\overline D:
 sx^n\in P\text{ for some }s\in S, n\geq1\}.
\]
More precisely, whenever $sx^n\in D$ with $s\in S$ and $n\geq1$,
one has $x\in Q(P)$ if and only if $sx^n\in P$.
\end{theorem}
\begin{proof}
Put $A=D_S$ and $B=(\overline D)_S$.
By Proposition~\ref{p:localized-root-extension}, $A\subseteq B$
is an integral root extension. Lying over makes the map
$\Spec(B)\longrightarrow\Spec(A)$ surjective. If $\mathfrak q$
lies over $\mathfrak p$ and $b\in B$, choose $n\geq1$ with
$b^n\in A$. Then $b\in\mathfrak q$ holds precisely when
$b^n\in\mathfrak p$. Thus there is at most one prime over each
$\mathfrak p$. An integral map on spectra is closed, so this
continuous bijection is a homeomorphism.

The standard localization homeomorphisms identify these spectra
with $U_S(\overline D)$ and $U_S(D)$. If $Q$ corresponds to $P$
and $sx^n\in D$, then $s\notin Q$ and primality gives
$x\in Q$ precisely when $sx^n\in Q\cap D=P$.
Proposition~\ref{p:localized-root-extension} supplies such a pair
$(s,n)$ for every $x\in\overline D$, proving the formula and
its independence of the choice of clearing data.
\end{proof}

Descent from the AGCD ring $(\overline D)_S$ to $\overline D$
still requires control of denominators. The next statement makes
that requirement explicit.

\begin{theorem}\label{t:integral-closure-descent}
Let $D$ be a domain, let $S\subseteq D^\ast$ be multiplicative,
and assume that $D_S$ is AGCD. Suppose that
\[
 \text{for each }0\ne x\in\overline D, \text{some }s\in S
 \text{ satisfies }
 s\bigl(x(\overline D)_S\cap\overline D\bigr)
 \subseteq x\overline D.\tag{$C_{\overline D}$}
\]
Then $\overline D$ is almost $S$-GCD. In particular, the conclusion
holds for an almost $S$-GCD domain satisfying this additional
condition on its normalization.
\end{theorem}
\begin{proof}
Proposition~\ref{p:localized-root-extension} shows that
$(\overline D)_S$ is AGCD. Apply Theorem~\ref{t3} to
$(\overline D,S)$.
\end{proof}

\begin{corollary}\label{c:integral-closure-tfinite}
Let $D$ be a domain and let $S\subseteq D^\ast$ be multiplicative.
Assume that $D_S$ is AGCD. If every $t$-ideal of $\overline D$
is of $v$-finite type, in particular if $\overline D$ is Mori,
then $\overline D$ is almost $S$-GCD.
\end{corollary}
\begin{proof}
The contraction $x(\overline D)_S\cap\overline D$ is a $t$-ideal.
A finite $v$-basis for this contraction allows its finitely many
denominators to be cleared by one element of $S$. Taking
$v$-closures then gives $(C_{\overline D})$.
Apply Theorem~\ref{t:integral-closure-descent}.
\end{proof}

\begin{corollary}\label{c:integral-closure-splitting}
Let $D$ be a domain and let $S\subseteq D^\ast$ be multiplicative.
Assume that $D_S$ is AGCD and that $S$, considered in
$\overline D$, is a splitting set. Then $\overline D$ is
almost $S$-GCD.
\end{corollary}
\begin{proof}
For a splitting set, each $x\in\overline D\setminus\{0\}$ can
be written $x=st$ with $s\in S$ and $t$ $v$-coprime to $S$.
Then $x(\overline D)_S\cap\overline D=t\overline D$.
Multiplication by $s$ proves $(C_{\overline D})$, so the result
follows from Theorem~\ref{t:integral-closure-descent}.
\end{proof}

\subsection{Content and polynomial descent}

For a polynomial $F$ over $K$ in any set of indeterminates,
$c(F)$ denotes the fractional ideal of $D$ generated by its
coefficients. The following standard content facts, together
with a pure-power observation, will be used explicitly.

\begin{lemma}\label{l:normal-content}
Let $D$ be an integrally closed domain with quotient field $K$.
\begin{enumerate}
\item For nonzero polynomials $F,G$ over $K$ in any set of
indeterminates, $c(FG)_v=(c(F)c(G))_v$.
\item If $A=(a_0,\ldots,a_r)$ is a nonzero finitely generated
fractional ideal and $n\geq1$, then
$(A^n)_v=(a_0^n,\ldots,a_r^n)_v$.
\item If $0\ne H\in D[X]$ and $0\ne s\in D$, then
$(H,s)_{v_{D[X]}}=D[X]$ if and only if $(c(H)+sD)_v=D$.
\end{enumerate}
\end{lemma}
\begin{proof}
(1) \cite[Proposition 34.8]{Gilmer}.
%An integrally closed domain is the intersection of its valuation
%overrings in $K$ \cite{Gilmer}. Over a valuation domain $V$, the
%content of a polynomial is principal and the usual Gauss argument
%gives $c_V(FG)=c_V(F)c_V(G)$. If $c(FG)\subseteq xD$ for
%$x\in K^\ast$, it follows that $c(F)c(G)\subseteq xV$ for every
%such $V$, and hence $c(F)c(G)\subseteq xD$.
%Intersecting over all principal fractional ideals containing
%$c(FG)$ proves $c(F)c(G)\subseteq c(FG)_v$; the reverse
%inclusion after taking $v$-closures is immediate.
%For several variables, alternatively apply the one-variable
%formula after a Kronecker substitution whose base exceeds
%all the multidegrees of $FG$. Only finitely many variables
%occur, so this also covers an arbitrary variable set.

(2) Put $I=(a_0^n,\ldots,a_r^n)$. If $I\subseteq xD$ and
$u=\prod_i a_i^{e_i}$ with $\sum_i e_i=n$, then
$(u/x)^n=\prod_i(a_i^n/x)^{e_i}\in D$.
Since $D$ is root closed, $u\in xD$. Intersecting all such
principal fractional ideals gives $A^n\subseteq I_v$.
Since $I\subseteq A^n$, their $v$-closures agree.

(3) Set $L=c(H)+sD$. If $z\in(H,s)^{-1}$, then $sz\in D[X]$
forces $z\in K[X]$. By (1), the inclusion $Hz\in D[X]$
implies $c(H)c(z)\subseteq D$. Together with $sc(z)\subseteq D$,
this says that every coefficient of $z$ lies in $L^{-1}$.
The converse follows by multiplying coefficients. Thus
$(H,s)^{-1}=L^{-1}[X]$, where the right side denotes the
polynomials with coefficients in $L^{-1}$.
Consequently, $(H,s)^{-1}=D[X]$ precisely when $L^{-1}=D$,
which proves the assertion.
\end{proof}

Descent to constants needs no normality assumption.

\begin{proposition}\label{p:poly-descent-to-D}
Let $D$ be a domain, let $S\subseteq D^\ast$ be multiplicative,
and let $\mathcal X$ be any set of indeterminates.
If $D[\mathcal X]$ is almost $S$-GCD, then $D$ is almost $S$-GCD.
\end{proposition}
\begin{proof}
Put $R=D[\mathcal X]$. For $a,b\in D^\ast$, choose $n\geq1$,
$s\in S$, and $0\ne F\in R$ such that
$s(a^nR\cap b^nR)\subseteq FR\subseteq a^nR\cap b^nR$.
The intersection equals $(a^nD\cap b^nD)R$ and contains a
nonzero constant. Thus $FR$ contains a nonzero constant, which
forces $F\in D$ by degree considerations. Contracting the
inclusions to $D$ proves the assertion by Theorem~\ref{t1}.
\end{proof}

The general polynomial criterion separates the denominator condition
from the classical root-extension obstruction.

\begin{theorem}\label{t:poly-localization-criterion}
Let $D$ be a domain, let $S\subseteq D^\ast$ be multiplicative,
and put $R=D[X]$. Assume that for every $0\ne f\in R$ there
is $s\in S$ such that $s(fR_S\cap R)\subseteq fR$.
The following conditions are equivalent:
\begin{enumerate}
\item $D[X]$ is almost $S$-GCD;
\item $D_S[X]$ is AGCD;
\item $D_S$ is AGCD and $D_S[X]\subseteq\overline{D_S}[X]$
is a root extension.
\end{enumerate}
\end{theorem}
\begin{proof}
Since $R_S=D_S[X]$, Theorem~\ref{t3} gives the equivalence of
(1) and (2). The equivalence of (2) and (3) is the classical
polynomial theorem for AGCD domains \cite[Theorem~3.4]{ATZ}.
\end{proof}

\begin{corollary}\label{c:poly-t-finite}
Let $D$ be a domain and let $S\subseteq D^\ast$ be multiplicative.
If every $t$-ideal of $D[X]$ is of $v$-finite type, in particular
if $D[X]$ is Mori, then the three conditions in
Theorem~\ref{t:poly-localization-criterion} are equivalent.
\end{corollary}
\begin{proof}
Each contraction $fD_S[X]\cap D[X]$ is a $t$-ideal.
Clearing denominators on a finite $v$-basis gives the hypothesis
of Theorem~\ref{t:poly-localization-criterion}.
\end{proof}

\begin{corollary}\label{c:poly-integrally-closed-localization}
Under the contraction hypothesis of
Theorem~\ref{t:poly-localization-criterion}, if $D_S$ is
integrally closed, then $D[X]$ is almost $S$-GCD if and only
if $D_S$ is AGCD.
\end{corollary}
\begin{proof}
The root extension in condition (3) is then the identity.
\end{proof}

\begin{corollary}\label{c:necessary-poly-conditions}
Let $D$ be a domain and let $S\subseteq D^\ast$ be multiplicative.
If $D[X]$ is almost $S$-GCD, then $D$ is almost $S$-GCD,
$D_S$ is AGCD, and $D_S[X]\subseteq\overline{D_S}[X]$ is a
root extension. No contraction hypothesis is needed for these
necessary conditions.
\end{corollary}
\begin{proof}
Apply Proposition~\ref{p:poly-descent-to-D}, localize at $S$,
and use \cite[Theorem~3.4]{ATZ}.
\end{proof}

Over a normal coefficient domain, the contraction hypothesis for
polynomials is equivalent to its constant version.

\begin{proposition}\label{p:normal-polynomial-contraction}
Let $D$ be an integrally closed domain, let $S\subseteq D^\ast$
be multiplicative, and let $\mathcal X$ be any set of
indeterminates. Put $R=D[\mathcal X]$. The following are equivalent:
\begin{enumerate}
\item for every $0\ne d\in D$, some $s\in S$ satisfies
$s(dD_S\cap D)\subseteq dD$;
\item for every $0\ne f\in R$, some $s\in S$ satisfies
$s(fR_S\cap R)\subseteq fR$.
\end{enumerate}
\end{proposition}
\begin{proof}
Condition (2) implies (1) by taking constants.
Assume (1), fix $0\ne f\in R$, and choose a nonzero coefficient
$a$ of $f$. Choose $s\in S$ with
$s(aD_S\cap D)\subseteq aD$.
For $h\in fR_S\cap R$, write $h=fg$ with
$g\in D_S[\mathcal X]$.
Lemma~\ref{l:normal-content} gives
$c(f)c(g)\subseteq c(fg)_v\subseteq D$.
Every coefficient $b$ of $g$ therefore satisfies
$ab\in aD_S\cap D$, so $sab\in aD$ and $sb\in D$.
Hence $sg\in R$ and $sh\in fR$.
The chosen $s$ depends only on $f$, so it works also when
$\mathcal X$ is infinite.
\end{proof}

\begin{corollary}\label{c:normal-polynomial-arbitrary}
\label{c:poly-transfer-two-descent}
Let $D$ be an integrally closed domain and let
$S\subseteq D^\ast$ satisfy condition (1) of
Proposition~\ref{p:normal-polynomial-contraction}.
For every set of indeterminates $\mathcal X$, the domain
$D[\mathcal X]$ is almost $S$-GCD if and only if $D$ is
almost $S$-GCD.
\end{corollary}
\begin{proof}
The proposition and Theorem~\ref{t3} reduce the assertions to
$D_S[\mathcal X]$ and $D_S$ being AGCD, respectively.
Since $D_S$ is integrally closed, the classical polynomial
criterion \cite[Theorem~3.4]{ATZ} gives their equivalence for
one variable, and iteration handles finitely many variables.
Any pair of polynomials in an arbitrary variable set lies in
a finite polynomial subring. A principal intersection there
remains principal after adjoining the remaining variables,
since a polynomial extension is flat. This proves the result
for arbitrary $\mathcal X$ as well.
\end{proof}

For PvMDs, the class-group criterion removes the contraction
hypothesis altogether.

\begin{corollary}\label{c:pvmd-polynomial-arbitrary}
Let $D$ be a PvMD, let $S\subseteq D^\ast$ be multiplicative,
and let $\mathcal X$ be any set of indeterminates. Then
$D[\mathcal X]$ is almost $S$-GCD if and only if $D$ is
almost $S$-GCD.
\end{corollary}
\begin{proof}
Polynomial rings in finitely many variables over a PvMD are
PvMDs, and a PvMD is integrally closed \cite{Kang}.
Theorem~\ref{t:apvmd-relative-characterization}, applied to
$D$ and these polynomial rings, reduces the assertion to
polynomial stability of the integrally closed AGCD domain
$D_S$. For an arbitrary variable set, use
Proposition~\ref{p:poly-descent-to-D} for descent and extend
intersections from finite polynomial subrings for ascent,
as in the preceding proof.
\end{proof}

\subsection{Almost splitting sets}

The splitting conditions below originated in \cite{ATZ}.
For a subset $S\subseteq D^\ast$, set
$S^{\perp_D}=\{t\in D^\ast:(t,s)_v=D\text{ for all }s\in S\}$.
A saturated multiplicative set $S$ is \emph{almost splitting}
if each $d\in D^\ast$ admits $d^n=st$ with $n\geq1$,
$s\in S$, and $t\in S^{\perp_D}$.
It is \emph{almost lcm-splitting} if, in addition, for every
$d\in D^\ast$ and $s\in S$, some $d^nD\cap s^nD$ is
principal. The sets $S$ and $S^{\perp_D}$ are multiplicative;
for the latter assertion, use the primality of maximal
$t$-ideals.

The finite form of the almost lcm condition needs an explicit
synchronization of exponents.

\begin{lemma}\label{l:finite-lcm-form}
Let $S$ be an almost lcm-splitting set of a domain $D$.
For every finite family $s_0,\ldots,s_r\in S$, there are
$m\geq1$ and $s\in S$ with $(s_0^m,\ldots,s_r^m)_v=sD$.
\end{lemma}
\begin{proof}
If $(u_1,\ldots,u_k)_v=D$, then
$(u_1^q,\ldots,u_k^q)_v=D$ for every $q\geq1$:
otherwise a maximal $t$-ideal would contain all the $u_i^q$
and hence all the $u_i$. Thus an equality
$(x_1^m,\ldots,x_k^m)_v=dD$ persists after raising its
normalized generators and $d$ to a common positive power.

Proceed by induction on the number of elements. For two
elements, the almost lcm condition and the equivalence between
principal intersections and principal divisorial sums give
the assertion. Its generator divides an element of $S$ and
therefore belongs to $S$ by saturation. If
$(s_0^q,\ldots,s_{r-1}^q)_v=dD$ with $d\in S$, apply the
two-element case to $d,s_r^q$ and choose $k\geq1$ with
$(d^k,s_r^{qk})_v=sD$, $s\in S$.
The initial observation gives
$(s_0^{qk},\ldots,s_r^{qk})_v=sD$.
\end{proof}

We now relate almost lcm-splitting to polynomial almost splitting.
The proof uses normal content identities and the primality of
maximal $t$-ideals.

\begin{theorem}\label{t:almost-lcm-polynomial-splitting}
Let $D$ be an integrally closed domain and let
$S\subseteq D^\ast$ be saturated and multiplicative.
Then $S$ is almost lcm-splitting in $D$ if and only if $S$
is almost splitting in $D[X]$.
\end{theorem}
\begin{proof}
Put $T=S^{\perp_D}$. Suppose first that $S$ is almost splitting
in $D[X]$. Applying its definition to constants and using
degrees shows that $S$ is almost splitting in $D$.
Given $d\in D^\ast$ and $s_0\in S$, write
$(d+s_0X)^n=sH$, where $n\geq1$, $s\in S$, and $H$ is
$v$-coprime to every element of $S$.
The leading coefficient $s_0^n/s$ of $H$ belongs to $D$, and
then to $S$ by saturation. Apply Lemma~\ref{l:normal-content}(3) to $H$ and this leading
coefficient, which belongs to $S\cap c(H)$. It gives $c(H)_v=D$. Parts (1) and (2) of that lemma now give
\[
 (d^n,s_0^n)_v=c((d+s_0X)^n)_v=sD.
\]
Thus $d^nD\cap s_0^nD$ is principal, proving the almost
lcm condition.

Conversely, suppose that $S$ is almost lcm-splitting and fix
$0\ne f\in D[X]$. List its nonzero coefficients as
$a_0,\ldots,a_r$ and put $A=c(f)$.
Almost splitting and a common multiple of the resulting
exponents give $a_i^q=s_it_i$ with $s_i\in S$ and $t_i\in T$.
Choose $m\geq1$ and $s\in S$ such that
$(s_0^m,\ldots,s_r^m)_v=sD$, by Lemma~\ref{l:finite-lcm-form}.
Set $u_i=s_i^m/s\in D$, $N=qm$, and
$B=(u_0t_0^m,\ldots,u_rt_r^m)$.
Then $(u_0,\ldots,u_r)_v=D$ and
$(a_0^N,\ldots,a_r^N)=sB$.

For each $s'\in S$, we have $(B+s'D)_v=D$.
Indeed, otherwise a maximal $t$-ideal $M$ would contain
$B+s'D$. The equality $(t_i,s')_v=D$ forces $t_i\notin M$.
As $u_it_i^m\in M$, primality then gives every $u_i\in M$,
contradicting $(u_0,\ldots,u_r)_v=D$.
By Lemma~\ref{l:normal-content}(2),
$(A^N)_v=(a_0^N,\ldots,a_r^N)_v=sB_v\subseteq sD$.
Hence $f^N=sH$ with $H\in D[X]$, and part (1) gives
\[
 c(H)_v=s^{-1}c(f^N)_v=s^{-1}(A^N)_v=B_v.
\]
Thus $(c(H)+s'D)_v=D$ for every $s'\in S$.
Part (3) shows that $H$ is $v$-coprime to $S$ in $D[X]$,
as required.
\end{proof}

The next example verifies directly the need for normality in
this equivalence.

\begin{example}\label{e:splitting-needs-normality}
Let $k=\mathbb F_2$, let $L=\mathbb F_4$, choose
$\alpha\in L\setminus k$, and put $B=L[t]$, $D=k+tB$,
and $S=D^\ast$. Then $D$ is AGCD, $S$ is almost lcm-splitting
in $D$, and $S$ is not almost splitting in $D[X]$.

Every $h\in B$ satisfies $h^3\in D$, since $h(0)^3\in k$.
Given $a,b\in D^\ast$, choose a gcd $d$ in the PID $B$ and
write $a=du$, $b=dv$ with $(u,v)B=B$.
If $x\in K=L(t)$ satisfies $xu^3,xv^3\in D$, then
$(u^3,v^3)B=B$ gives $x\in B$.
At least one of $u(0),v(0)$ is nonzero, so at least one of
$u(0)^3,v(0)^3$ equals $1$. It follows that $x(0)\in k$,
and hence $x\in D$.
Thus $(u^3,v^3)^{-1}=D$ and
$(a^3,b^3)_v=d^3D$. In particular, $D$ is AGCD with the
uniform exponent $3$. The set $S$ is trivially almost
splitting in $D$, and the AGCD property supplies its almost
lcm condition.

Put $F=X+\alpha$ and $f=tF\in D[X]$.
Suppose that $f^n=sH$ with $n\geq1$, $s\in S$, and
$H\in D[X]$ $v$-coprime to every element of $S$.
Write $H=cF^n$, where $c=t^n/s\in L(t)$.
The leading coefficient gives $c\in D$.
Writing $n=2^em$ with $m$ odd, the coefficient of
$X^{n-2^e}$ in $F^n$ is $\alpha^{2^e}\notin k$.
Since $c\alpha^{2^e}\in D$, evaluation at $t=0$ forces
$c(0)=0$, so $c\in tB$.
Consequently, $\alpha H\in D[X]$ and $\alpha t\in D$,
whereas $\alpha\notin D[X]$.
Thus $\alpha\in(H,t)^{-1}\setminus D[X]$, and
$(H,t)_{v_{D[X]}}\ne D[X]$, a contradiction.
\end{example}

\begin{remark}\label{r:integrally-closed-needed-q7}
The domain in Example~\ref{e:splitting-needs-normality} is not
integrally closed: $\alpha\in\qf(D)$ is integral over $D$ but
does not belong to $D$.
Thus the implication from almost lcm-splitting in $D$ to
almost splitting in $D[X]$ fails without the normality
hypothesis. This is consistent with the distinction between
AGCD domains and polynomial almost-splitting conditions
studied in \cite{AndersonChang2007}.
\end{remark}

We conclude with the consequence of the classical descent theorem
for almost lcm-splitting sets.

\begin{remark}\label{r:root-closed-lcm}
Let $D$ be root closed and let $S$ be an almost lcm-splitting
set of $D$. Then $D$ is almost $S$-GCD if and only if it is AGCD.
Indeed, almost $S$-GCD implies that $D_S$ is AGCD; the reverse
descent to $D$ follows from \cite[Corollary~2.10]{ATZ}.
The converse implication is immediate from the definition.
Whenever these conditions hold, $D$ is also integrally closed,
since the normalization of an AGCD domain is a root extension
of that domain.
\end{remark}

\begin{theorem}\label{t:poly-almost-S-GCD-lcm}
Let $D$ be a root-closed domain and let $S$ be an almost
lcm-splitting set of $D$. For any set of indeterminates
$\mathcal X$, the domain $D[\mathcal X]$ is almost $S$-GCD
if and only if $D$ is almost $S$-GCD.
\end{theorem}
\begin{proof}
Descent follows from Proposition~\ref{p:poly-descent-to-D}.
For ascent, Remark~\ref{r:root-closed-lcm} shows that $D$
is an integrally closed AGCD domain. Its polynomial rings
are AGCD by \cite[Theorem~3.4]{ATZ}, using iteration and
flat extension from finite polynomial subrings for arbitrary
$\mathcal X$. They are therefore almost $S$-GCD.
\end{proof}

\section{Conductor transfer and pullbacks}\label{sec:conductor}

A common conductor supplies the uniform multiplier needed to transfer
almost $S$-GCD in both directions. This yields pullback and
normalization results without a separate principal-contraction
hypothesis.
For domains $D\subseteq T\subseteq K$, the \emph{conductor} is
$(D:T)=\{c\in D:cT\subseteq D\}$. It is an ideal of both rings.

\begin{theorem}\label{t:conductor-transfer}
Let $D\subseteq T$ be domains with the same quotient field, let
$S\subseteq D^\ast$ be multiplicative, and suppose that
$S\cap(D:T)\neq\varnothing$. Then $D$ is almost $S$-GCD if and only
if $T$ is almost $S$-GCD, with $S$ viewed as a subset of $T$.
\end{theorem}

\begin{proof}
Fix $c\in S\cap(D:T)$. First suppose that $T$ is almost $S$-GCD.
For $a,b\in D^\ast$, choose $n\geq1$, $s\in S$, and $d\in T^\ast$
with $sJ_T\subseteq dT\subseteq J_T$, where
$J_T=a^nT\cap b^nT$. Set $J_D=a^nD\cap b^nD$.
Then $J_D\subseteq J_T$ and $cJ_T\subseteq J_D$, so
$c^2sJ_D\subseteq cdD\subseteq J_D$.

Conversely, suppose that $D$ is almost $S$-GCD and take
$a,b\in T^\ast$. Apply the defining property to $ca,cb\in D$ to
find $n\geq1$, $s\in S$, and $d\in D^\ast$ such that
$sJ_D\subseteq dD\subseteq J_D$, where
$J_D=(ca)^nD\cap(cb)^nD$.
Put $J_T=(ca)^nT\cap(cb)^nT$. Since $cJ_T\subseteq J_D$,
$csJ_T\subseteq dT\subseteq J_T$.
Dividing this sandwich by $c^n$ shows that
$a^nT\cap b^nT$ is $S$-principal.
\end{proof}

\begin{remark}\label{r:conductor-S-GCD}
The same conductor-transfer theorem holds for $S$-GCD domains.
Indeed, in its proof all exponents can be taken equal to $1$.
The two-element $S$-principality condition for intersections
characterizes $S$-GCD domains \cite[Theorem~2.1]{AHZ}.
\end{remark}

We state the pullback application with all maps and conductor
hypotheses explicit.

\begin{corollary}\label{p:milnor-meets-M}
Let $T$ and $F$ be domains, let $\pi:T\to F$ be a surjective
homomorphism with $M=\ker(\pi)\ne0$, and let $A$ be a subdomain
of $F$. Form the pullback $R=\pi^{-1}(A)$:
\[
 \xymatrix{
 R \ar[r] \ar[d] & T \ar[d]^{\pi}\\
 A \ar@{^{(}->}[r] & F .
 }
\]
Let $S\subseteq R^\ast$ be multiplicative with $S\cap M\ne\varnothing$.
Then $R_S=T_S$, and $R$ is almost $S$-GCD if and only if $T$ is
almost $S$-GCD.
If $T$ is an APvMD, these conditions are equivalent to $T_S$
being AGCD. The latter equivalence also holds if, instead,
$(C_{R,S})$ holds.
\end{corollary}

\begin{proof}
The ideal $M$ is contained in $(R:T)$.
For $c\in S\cap M$ and $t\in T$, one has $ct\in M\subseteq R$,
so $t=(ct)/c\in R_S$. Thus $R_S=T_S$; the same argument also
shows $\qf(R)=\qf(T)$.
Apply Theorem~\ref{t:conductor-transfer}. If $T$ is an APvMD,
apply Theorem~\ref{t:apvmd-relative-characterization} to $T$.
If $(C_{R,S})$ holds, apply Theorem~\ref{t3} to $R$.
\end{proof}

\begin{corollary}\label{c:normalization-conductor}
Let $D$ be a domain, $S\subseteq D^\ast$ multiplicative, and
suppose $S\cap(D:\overline D)\ne\varnothing$.
Then $D$ is almost $S$-GCD if and only if $\overline D$ is almost
$S$-GCD. If $\overline D$ is a PvMD, these conditions are also
equivalent to $D_S$ being AGCD.
\end{corollary}

\begin{proof}
Apply conductor transfer to $D\subseteq\overline D$.
A conductor element in $S$ gives
$D_S=(\overline D)_S$, and the last assertion follows from
Theorem~\ref{t:apvmd-relative-characterization}.
\end{proof}

The conductor theorem yields nonnormal polynomial examples directly.

\begin{example}\label{e:cusp-conductor}
Let $D=\QQ[t^2,t^3]$ and $S=\langle t^2\rangle$.
Since $t^2\QQ[t]\subseteq D$, the conductor meets $S$.
For any set of indeterminates $\mathcal X$, the same element
lies in $(D[\mathcal X]:\QQ[t][\mathcal X])$.
The upper polynomial ring is GCD, so
$D[\mathcal X]$ is $S$-GCD by
Remark~\ref{r:conductor-S-GCD}.
Nevertheless $D$ is not AGCD:
$1+t$ is integral over $D$, while the coefficient of $t$ in
$(1+t)^n$ is $n\ne0$ for every $n\geq1$. Thus
$D\subseteq\overline D=\QQ[t]$ is not a root extension.
\end{example}

\begin{example}\label{e:field-pullback}
For any extension of fields $k\subseteq L$, the domain
$D=k+XL[X]$ is $\langle X\rangle$-GCD.
Indeed, $XL[X]$ is a common conductor ideal and $L[X]$ is a PID.
This statement does not require $L/k$ to be algebraic or a root
extension.
\end{example}

\section{Cyclic localizations and finite principal covers}\label{sec:sublocal}

The sublocal condition asks for every nonunit denominator in $S$
to work separately. We first record its elementary behavior and a
strict separation from the global condition. We then connect finite
principal covers with invertible power intersections and the Picard
group.

\subsection{The sublocal condition}

\begin{definition}\label{d:sublocal-almost-S-GCD}
Let $D$ be a domain and $S\subseteq D^\ast$ multiplicative.
Put $S^\circ=S\setminus U(D)$. The domain $D$ is
\emph{sublocally almost $S$-GCD} if it is almost
$\langle s\rangle$-GCD for every $s\in S^\circ$.
\end{definition}

The exclusion of units is part of the definition. Requiring the
condition for every $s\in S$, including $1$, would make it
equivalent to AGCD. On the other hand, if $S\subseteq U(D)$,
Definition~\ref{d:sublocal-almost-S-GCD} is vacuous and holds for
every domain. Thus the implication to the global almost
$S$-GCD condition requires $S^\circ\ne\varnothing$.

The cyclic condition is controlled by divisibility up to powers.

\begin{lemma}\label{l:cyclic-power-divisibility}
Let $D$ be a domain and let $s,t\in D^\ast$.
For every $m\geq1$, the almost $\langle s\rangle$-GCD and
almost $\langle s^m\rangle$-GCD conditions are equivalent.
If $s$ divides a positive power of $t$, the almost
$\langle s\rangle$-GCD condition implies the almost
$\langle t\rangle$-GCD condition.
In particular, the two cyclic conditions agree if
$\sqrt{sD}=\sqrt{tD}$.
\end{lemma}

\begin{proof}
If $s^kI\subseteq dD\subseteq I$, replace $k$ by a larger
multiple of $m$ to obtain an $\langle s^m\rangle$-witness.
The reverse implication follows by inclusion of multiplicative
sets. If $t^q=su$, then
$t^{qk}I=u^ks^kI\subseteq dD$, so the same principal ideal
works for $\langle t\rangle$. Equality of the radicals gives
divisibility by a positive power in both directions.
\end{proof}

Several permanence properties can now be stated together.

\begin{proposition}\label{p:sublocal-properties}
\label{c:sublocal-generators}\label{c:principal-sublocal}
\label{p:sublocal-implies-global}\label{p:AGCD-implies-sublocal}
\label{p:sublocal-localization}\label{p:sublocal-polynomial-descent}
Let $D$ be a domain and $S,U\subseteq D^\ast$ multiplicative.
\begin{enumerate}
\item Every AGCD domain is sublocally almost $S$-GCD.
If $S^\circ\ne\varnothing$, every sublocally almost $S$-GCD
domain is almost $S$-GCD.
\item If $S$ is generated by $\Gamma\subseteq D^\ast$, then
$D$ is sublocally almost $S$-GCD if and only if it is almost
$\langle\gamma\rangle$-GCD for every nonunit $\gamma\in\Gamma$.
For a nonzero nonunit $s$, the sublocal and global conditions
for $S=\langle s\rangle$ therefore coincide.
\item If $D$ is sublocally almost $S$-GCD, then $D_U$ is
sublocally almost $S_U$-GCD, where $S_U=\{s/1:s\in S\}$.
\item If $D[\mathcal X]$ is sublocally almost $S$-GCD, with $S$
viewed as constants, then $D$ is sublocally almost $S$-GCD.
\end{enumerate}
\end{proposition}

\begin{proof}
The first assertion follows from
Proposition~\ref{p:permanence}; for the second part of (1),
choose any $s\in S^\circ$ and use $\langle s\rangle\subseteq S$.
For (2), every nonunit finite product of generators has a
nonunit factor in $\Gamma$. Apply
Lemma~\ref{l:cyclic-power-divisibility}.
For (3), if $s/1$ is a nonunit in $D_U$, then $s$ was a nonunit
in $D$; localize the almost $\langle s\rangle$-GCD condition.
For (4), nonunit constants remain nonunits in a polynomial
ring, and Proposition~\ref{p:poly-descent-to-D} applies to
each cyclic set.
\end{proof}

The localization criteria of the preceding sections give
corresponding sublocal criteria.

\begin{proposition}\label{p:cyclic-localization-criterion}
\label{c:sublocal-tfinite}
Let $D$ be a domain and $S\subseteq D^\ast$ multiplicative.
Assume that at least one of the following holds:
\begin{enumerate}
\item $(C_{D,\langle s\rangle})$ holds for every $s\in S^\circ$;
\item $D$ is Mori;
\item every $aD\cap bD$, $a,b\in D^\ast$, is of $v$-finite type;
\item $D$ is an APvMD.
\end{enumerate}
Then $D$ is sublocally almost $S$-GCD if and only if $D_s$ is
AGCD for every $s\in S^\circ$.
\end{proposition}

\begin{proof}
Apply, respectively, Theorem~\ref{t3},
Corollary~\ref{c:tfinite-localization-criterion},
Corollary~\ref{c:finite-conductor}, or
Theorem~\ref{t:apvmd-relative-characterization}
to each cyclic multiplicative set.
\end{proof}

For Krull domains the cyclic obstruction is read directly from
prime-divisor classes.

\begin{corollary}\label{c:krull-sublocal-supports}
Let $D$ be a Krull domain and $S\subseteq D^\ast$ multiplicative.
Then $D$ is sublocally almost $S$-GCD if and only if, for every
$s\in S^\circ$, the classes of the height-one prime divisors of
$s$ span $\Cl(D)\otimes_{\ZZ}\QQ$.
\end{corollary}

\begin{proof}
Apply Theorem~\ref{t:krull-boundary-criterion} to
$\langle s\rangle$. A height-one prime meets this cyclic set
exactly when it contains $s$.
\end{proof}

The global condition can be strictly weaker, although it agrees
with the sublocal condition for a cyclic set generated by a nonunit.

\begin{example}\label{e:global-not-sublocal}
Let $D=\QQ[t^2,t^3,Y]$ and let
$S=\{t^{2a}Y^b:a,b\geq0\}$.
Then $D$ is Noetherian and
$D_S=\QQ[t,t^{-1},Y,Y^{-1}]$ is a UFD, so $D$ is almost
$S$-GCD by Corollary~\ref{c:tfinite-localization-criterion}.
However, $D_Y=\QQ[t^2,t^3,Y,Y^{-1}]$ is not AGCD.
The element $1+t$ is integral over this ring, satisfying
$(Z-1)^2-t^2=0$, but no positive power belongs to it: the
coefficient of $t$ in $(1+t)^n$ is the nonzero rational number
$n$. The normalization therefore fails to be a root extension.
Since $Y\in S^\circ$, the domain $D$ is not sublocally almost
$S$-GCD.
\end{example}

Conductor transfer also gives a sublocal version, with nonunits
always measured in the lower ring.

\begin{proposition}\label{p:sublocal-milnor-radical}
Let $R\subseteq T\subseteq\qf(R)$ be domains and let
$S\subseteq R^\ast$ be multiplicative. Put
$S^\circ=S\setminus U(R)$ and assume
$S^\circ\subseteq\sqrt{(R:T)}$.
If $T$ is an APvMD, then $R$ is sublocally almost $S$-GCD
if and only if $T_s$ is AGCD for every $s\in S^\circ$.
The same conclusion holds if, instead of the APvMD hypothesis,
$(C_{R,\langle s\rangle})$ holds for every $s\in S^\circ$.
\end{proposition}

\begin{proof}
For each such $s$, some $s^m$ belongs to $(R:T)$.
Thus $R_s=T_s$, and conductor transfer for $\langle s\rangle$
identifies the almost cyclic conditions on $R$ and $T$.
If $T$ is an APvMD, apply
Theorem~\ref{t:apvmd-relative-characterization} to $T$.
Under the alternative contraction hypothesis apply
Theorem~\ref{t3} to $R$.
This also applies to the pullback in
Corollary~\ref{p:milnor-meets-M} whenever
$S^\circ\subseteq\sqrt M$.
\end{proof}

\subsection{Finite principal covers and the Picard group}

The \emph{Picard group} $\Pic(D)$ is the group of nonzero
invertible fractional ideals modulo principal fractional ideals.
A domain is an \emph{almost generalized GCD domain} (AGGCD domain)
if each pair $a,b\in D^\ast$ has some invertible intersection
$a^nD\cap b^nD$. This established class lies between AGCD domains
and the corresponding local conditions \cite{ACZlocal, Lewin}.
The next theorem explains how it arises from a finite principal
cover and why a class-group obstruction remains.

\begin{theorem}\label{t:principal-cover-picard}
Let $D$ be a domain and let $s_1,\ldots,s_r\in D^\ast$ satisfy
$(s_1,\ldots,s_r)D=D$. Suppose $D_{s_i}$ is AGCD for every $i$.
Then, for every $a,b\in D^\ast$, there is $n\geq1$ for which
$J=a^nD\cap b^nD$ is invertible and
\[
 a^{nm}D\cap b^{nm}D=J^m
 \qquad\text{for every }m\geq1.
\]
In particular, $D$ is AGGCD. Moreover, $D$ is AGCD if and only
if $\Pic(D)$ is torsion.
\end{theorem}

\begin{proof}
For each $i$, choose $n_i$ making
$a^{n_i}D_{s_i}\cap b^{n_i}D_{s_i}$ principal.
When a pairwise intersection is principal, its higher power
intersections are the corresponding ordinary powers: apply
Lemma~\ref{l:power-sandwich} with multiplier $1$ to the
associated divisorial sum, and then invert.
Thus a common multiple $n$ of the $n_i$ makes
$JD_{s_i}$ principal for all $i$.

Choose $d_i\in J$ generating $JD_{s_i}$ and put
$N=(d_1,\ldots,d_r)D\subseteq J$.
Then $(J/N)_{s_i}=0$ for every $i$.
For an element $x\in J/N$, choose powers $s_i^{k_i}$ that
annihilate $x$. These powers generate the unit ideal, since
the $s_i$ do. Hence $x=0$, so $J=N$ is finitely generated.
It is locally principal on the principal open cover
$\{\mathfrak p:s_i\notin\mathfrak p\}$, and is consequently
invertible.

For every $m$, the displayed power identity holds after
localizing at each $s_i$, by the same principal-intersection
calculation. Equality on this cover gives equality over $D$.
If $\Pic(D)$ is torsion, choose $m$ such that $J^m$ is
principal. The displayed identity then proves the AGCD
condition for $a,b$. Conversely, the $t$-class group of an
AGCD domain is torsion \cite[Theorem~3.4]{AZ-Bez}, and
$\Pic(D)$ embeds into $\Cl_t(D)$ because every invertible
fractional ideal is a $t$-invertible $t$-ideal. Thus AGCD
implies torsion of $\Pic(D)$.
\end{proof}

\begin{corollary}\label{c:cyclic-cover}
Let $D$ be a domain and let $s_1,\ldots,s_r\in D^\ast$ generate
the unit ideal. If $D$ is almost $\langle s_i\rangle$-GCD for
every $i$, then $D$ is AGGCD. If, in addition, $\Pic(D)$ is
torsion, then $D$ is AGCD.
In particular, if $D$ is nonlocal, $S=D^\ast$, and $D$ is
sublocally almost $S$-GCD, then $D$ is AGGCD; it is AGCD
whenever $\Pic(D)$ is torsion. In particular, the latter conclusion
holds if $D$ is semilocal and nonlocal.
\end{corollary}

\begin{proof}
Localizing gives the hypotheses of
Theorem~\ref{t:principal-cover-picard}.
For the last assertion choose distinct maximal ideals
$M,N$ and elements $s\in M$, $t\in N$ with $s+t=1$.
Both are nonzero nonunits, so the sublocal hypothesis applies
to the principal cover given by $s,t$. For the semilocal assertion,
use that every invertible ideal of a semilocal ring is principal
\cite{Gilmer}, so $\Pic(D)=0$.
\end{proof}
\section*{Declarations}
\noindent\textbf{Availability of data and material.}
No datasets were generated or analysed during the current study.

\smallskip
\noindent\textbf{Competing interests.}
The authors declare that they have no competing interests.

\smallskip
\noindent\textbf{Authors' contributions.}
All authors contributed equally to the conception, writing, and revision
of the manuscript. All authors read and approved the final version
of the manuscript.

\smallskip
\noindent\textbf{Funding.}
H. Kim was supported by Basic Science Research Program through the
National Research Foundation of Korea (NRF) funded by the Ministry
of Education (grant no.~RS-2026-25589739).

\end{document}